\documentclass{amsart}

\usepackage{amssymb}
\usepackage{amsmath}
\usepackage{graphicx}

\newcommand{\erase}[1]{}

\newtheorem{theorem}{Theorem}[section]
\newtheorem{lemma}[theorem]{Lemma}
\newtheorem{proposition}[theorem]{Proposition}
\newtheorem{corollary}[theorem]{Corollary}

\newtheorem{_algorithm}[theorem]{Algorithm}
\newenvironment{algorithm}{\begin{_algorithm}\rm}{\hfill \rule{3pt}{6pt}\end{_algorithm}}

\newtheorem{_remark}[theorem]{\it Remark}
\newenvironment{remark}{\begin{_remark}\rm}{\end{_remark}}

\newtheorem{_example}[theorem]{Example}
\newenvironment{example}{\begin{_example}\rm}{\end{_example}}

\numberwithin{equation}{section}
\numberwithin{table}{section}
\numberwithin{figure}{section}

\newcommand{\C}{\mathord{\mathbb C}}

\renewcommand{\P}{\mathord{\mathbb  P}}
\newcommand{\Q}{\mathord{\mathbb  Q}}
\newcommand{\R}{\mathord{\mathbb R}}
\newcommand{\Z}{\mathord{\mathbb Z}}

\newcommand{\EEE}{\mathord{\mathcal E}}
\newcommand{\FFF}{\mathord{\mathcal F}}

\newcommand{\LLL}{\mathord{\mathcal L}}
\newcommand{\MMM}{\mathord{\mathcal M}}

\newcommand{\PPP}{\mathord{\mathcal P}}
\newcommand{\RRR}{\mathord{\mathcal R}}
\newcommand{\SSS}{\mathord{\mathcal S}}
\newcommand{\VVV}{\mathord{\mathcal V}}

\newcommand{\maprightsp}[1]{\; \smash{\mathop{\; \longrightarrow \; }\limits\sp{#1}}\; }
\newcommand{\isom}{\mathbin{\,\raise -.6pt\rlap{$\to$}\raise 3.5pt \hbox{\hskip .3pt$\mathord{\sim}$}\,\;}}

\newcommand{\set}[2]{\{\; {#1} \; \mid \; {#2} \;  \}}
\newcommand{\shortset}[2]{\{ {#1} \,|\, {#2}   \}}
\newcommand{\gen}[1]{\langle {#1}  \rangle}

\newcommand{\tensor}{\otimes}
\newcommand{\sprime}{\sp\prime}
\newcommand{\spprime}{\sp{\prime\prime}}
\newcommand{\sperp}{\sp{\perp}}

\newcommand{\dual}{\sp{\vee}}
\newcommand{\inv}{\sp{-1}}

\newcommand{\Hom}{\mathord{\mathrm {Hom}}}
\newcommand{\OG}{\mathord{\mathrm {O}}}
\newcommand{\Aut}{\operatorname{\mathrm {Aut}}\nolimits}

\newcommand{\Sing}{\operatorname{\mathrm {Sing}}\nolimits}
\newcommand{\rank}{\operatorname{\mathrm {rank}}\nolimits}

\newcommand{\quand}{\quad\textrm{and}\quad}

\newcommand{\intfL}[1]{\langle #1\rangle_{L}}
\newcommand{\intfS}[1]{\langle #1\rangle_{S}}
\newcommand{\intfW}[1]{\langle #1\rangle_{W}}
\newcommand{\intfLL}[1]{\langle #1\rangle_{\Lambda}}
\newcommand{\intfN}[1]{\langle #1\rangle_{N}}
\newcommand{\intfH}[1]{\langle #1\rangle_{H}}
\newcommand{\Weyl}[1]{W(#1)}
\newcommand{\va}{\mathord{\textbf{\itshape a}}}
\newcommand{\ve}{\mathord{\textbf{\itshape e}}}
\newcommand{\vt}{\mathord{\textbf{\itshape t}}}
\newcommand{\vv}{\mathord{\textbf{\itshape v}}}
\newcommand{\vw}{\mathord{\textbf{\itshape w}}}
\newcommand{\vae}{\va_{\varepsilon}}
\newcommand{\Lpsm}{\Lambda_{p, \sigma}^-}

\newcommand{\mE}{E_8^{(-1)}}
\newcommand{\pmE}{E_8^{(-p)}}
\newcommand{\Up}{U^{(p)}}
\newcommand{\Hp}{H^{(-p)}}

\newcommand{\thep}{7919}
\newcommand{\thepten}{17389}

\begin{document}

\title[Supersingular $K3$ surfaces  and Salem polynomials]
{Automorphisms of supersingular $K3$ surfaces and Salem polynomials}

\author{Ichiro Shimada}
\address{
Department of Mathematics, 
Graduate School of Science, 
Hiroshima University,
1-3-1 Kagamiyama, 
Higashi-Hiroshima, 
739-8526 JAPAN
}
\email{shimada@math.sci.hiroshima-u.ac.jp}

\thanks{Partially supported by
JSPS Grants-in-Aid for Scientific Research (C) No.~25400042%  and
%JSPS Grants-in-Aid for Scientific Research (A) No.~23244002%
}

\dedicatory{Dedicated to Professor Tetsuji Shioda on the occasion of his 75th birthday}

\begin{abstract}
We present a method to generate many automorphisms of a supersingular $K3$ surface in odd characteristic.
As an application,
we show that, if $p$ is an odd prime less than or equal to $\thep$,
then every supersingular $K3$ surface in characteristic $p$ has an automorphism 
whose characteristic polynomial  on the N\'eron--Severi lattice 
is a Salem polynomial of degree $22$.
For a supersingular $K3$ surface with Artin invariant $10$,
the same holds for odd primes less than or equal to $\thepten$.
\end{abstract}

%14J   (1973-now) Surfaces and higher-dimensional varieties [For analytic theory, see 32Jxx]
%14J28   (1980-now) K3 surfaces and Enriques surfaces
%14J   (1973-now) Surfaces and higher-dimensional varieties [For analytic theory, see 32Jxx]
%14J50   (1980-now) Automorphisms of surfaces and higher-dimensional varieties
%14Q   (1991-now) Computational aspects in algebraic geometry [See also 12Y05, 13Pxx, 68W30]
%14Q10   (1991-now) Surfaces, hypersurfaces
%37B  (2000-now) Topological dynamics [See also 54H20]
%37B40   (2000-now) Topological entropy
\maketitle

% main text starts here

\section{Introduction}\label{sec:Intro}
An irreducible monic polynomial $\phi(t)\in \Z[t]$ of even degree $2d>0$
is called a \emph{Salem polynomial} if 
$\phi(t)$ is reciprocal,
$\phi(t)=0$ has two positive real roots, and
the other $2d-2$ complex roots are located on $\shortset{z\in \C}{|z|=1}\setminus \{\pm 1\}$.
\par
The notion of Salem polynomials plays an important role in the study 
of dynamics of automorphisms of algebraic varieties.
We have the following fundamental theorem due to McMullen~\cite{MR1896103}.
See also~\cite{MR3129006} and~\cite[Proposition 3.1]{1406.2761}. 
\begin{theorem}[\cite{MR1896103}]\label{thm:cycsSalem}
Let $g$ be an automorphism of an algebraic  $K3$ surface $X$
defined over an algebraically closed field.
Then the characteristic polynomial of the action of $g$ on  
the N\'eron--Severi lattice $S_X$ of $X$ is
a product of cyclotomic polynomials and at most one Salem polynomial
counting with multiplicities.
\end{theorem}
A $K3$ surface $X$ defined over an algebraically closed field $k$ of characteristic $p>0$ is said to be \emph{supersingular}
if the rank of its N\'eron--Severi lattice $S_X$ is $22$.
%The supersingular $K3$ surfaces defined over the field $k$ constitute  a $9$-dimensional moduli.
We say that an automorphism $g$ of a supersingular $K3$ surface $X$ is \emph{of irreducible Salem type}
if the characteristic polynomial of the action of $g$  on  $S_X$
is a Salem polynomial of degree $22$.
\par
The purpose of this note is to report the following theorems,
which are the results of  computer-aided experiments.
By a \emph{double plane involution} of a $K3$ surface $X$ in characteristic not equal to $2$,
we mean an automorphism of $X$ of order $2$
induced by the Galois transformation  of a generically finite morphism $X\to \P^2$ of degree $2$.
\begin{theorem}\label{thm:main}
Let $p$ be an odd prime less than or equal to $\thep$.
Then every supersingular $K3$ surface $X$ in characteristic $p$ has 
a sequence  of 
double plane involutions $\tau_1, \dots,  \tau_l$ of length at most  $22$
such that their product $\tau_1\cdots \tau_l$ is 
an automorphism 
of irreducible Salem type.
\end{theorem}
Let $X$ be a supersingular $K3$ surface in characteristic $p>0$, and 
let $S_X\dual$ denote  the \emph{dual lattice} $\Hom(S_X, \Z)$ of $S_X$,
into which $S_X$ is embedded as a submodule of finite index by the intersection form of $S_X$.
Artin~\cite{MR0371899} showed that the discriminant group $S_X\dual/S_X$ of $S_X$ is isomorphic to 
$(\Z/p\Z)^{2\sigma}$,
where $\sigma$ is a positive integer less than or equal to $10$.
This integer $\sigma$ is called the \emph{Artin invariant} of $X$.
By the result of Ogus~\cite{MR563467, MR717616},
the supersingular $K3$ surfaces of Artin invariant $\le \sigma$ 
defined over an algebraically closed field $k$ constitute  a moduli of dimension $\sigma-1$,
and 
a supersingular $K3$ surface $X(p)$ with Artin invariant $1$ is unique up to isomorphism.
\par
For supersingular $K3$ surfaces  with Artin invariant  $\sigma=10$ in characteristic  $p$ with $11\le p\le \thepten$, 
we found a class of sequences  of 
double plane involutions whose product is \emph{frequently} of irreducible Salem type.
(See Section~\ref{sec:sigma10} for the detail.)
Using this class, we obtain the following theorem:
\begin{theorem}\label{thm:main2}
Let $p$ be an odd prime less than or equal to $\thepten$.
Then every supersingular $K3$ surface $X$ in characteristic $p$ with Artin invariant $10$ has 
a sequence  of 
double plane involutions  of length at most  $22$
such that their product is 
an automorphism 
of irreducible Salem type.
\end{theorem}
\par
The interest of an automorphism of irreducible Salem type  
stems from the following observation due to Esnault and Oguiso~\cite{1406.2761, 1411.0769}:
\begin{theorem}[\cite{1406.2761, 1411.0769}]\label{thm:nonliftable}
Let $g$ be an automorphism of a supersingular $K3$ surface $X$.
If the characteristic polynomial of the action of $g$ on  $S_X$ is irreducible,
then the pair $(X, g)$ can never be lifted to characteristic $0$.
\end{theorem}
Hence we obtain the following corollary.
\begin{corollary}\label{cor:main}
Let $X$ be a supersingular $K3$ surface in odd characteristic $p$ with Artin invariant $\sigma$.
Suppose that $p\le \thep$ or {\rm (}$\sigma=10$ and $p\le \thepten${\rm )}.
Then $X$  has an automorphism $g$
such that the pair $(X, g)$ can never be lifted to characteristic $0$.
\end{corollary}
Recently, several authors have studied  the non-liftability of automorphisms of supersingular $K3$ surfaces
by means of Salem polynomials.
See~\cite{1307.0361, 1406.2761, 1411.0769, 1502.06923}.
In particular,
the existence of a non-liftable automorphism has been  established
for a supersingular $K3$ surface $X(p)$  in characteristic $p$
\emph{with Artin invariant $1$}, 
except for the cases  $p= 7$ and $13$.
\begin{remark}
In~\cite{JangLocal}, the existence of a non-liftable  automorphism of $X(p)$ 
was proved for $p$ large enough by another method.
\end{remark}
Our main theorems  not only  fill the remaining cases $X(7)$ and $X(13)$ for  supersingular $K3$ surfaces with Artin invariant $1$,
but also  suggest that  this result can be extended  to  supersingular $K3$ surfaces with arbitrary Artin invariant,
at least in odd characteristics.
There exists no theoretical significance in
the bounds $p\le \thep$ in Theorem~\ref{thm:main} and $p\le \thepten$ in Theorem~\ref{thm:main2}.
We merely stopped our computations at the 1000th prime ($p=\thep$)
and the 2000th prime ($p=\thepten$).
\par
The main tool of the proof of Theorems~\ref{thm:main} and~\ref{thm:main2}
is the structure theorem of the N\'eron--Severi lattices of supersingular $K3$ surfaces $X$
due to Rudakov and Shafarevich~\cite{MR633161},
which states that 
 the isomorphism class of the lattice $S_X$ is uniquely determined by 
$p$ and the Artin invariant $\sigma$ of $X$.
%We construct a  lattice $\Lpsm$ isomorphic to $S_X$ explicitly in Section~\ref{sec:RS}.
\par
Let $X$ be a supersingular $K3$ surface $X$ in odd characteristic.
In this paper,  
we present a method to generate many matrix representations on $S_X$  
of double plane involutions of $X$.
Composing  some of these involutions,
we obtain an  automorphism of irreducible Salem type.
In order to produce double plane involutions,
we have to find the nef cone in $S_X\tensor\R$.
%,
%and hence we have to fix a standard fundamental domain
%of the action of the Weyl group $W(\Lpsm)$ on $\Lpsm\tensor\R$.
For this purpose,
we introduce  a notion of an \emph{ample list of vectors}.
(See Section~\ref{sec:lattice}  for the definitions.)
\par
The results of the experiments are presented in the author's web page~\cite{compdataIrredSalem}.
\par
Thanks are due to Professors Junmyeong Jang, Toshiyuki Katsura, Jonghae Keum,   Keiji Oguiso,  Matthias Sch\"utt
and Hirokazu Yanagihara 
for  stimulating discussions.
\section{Lattices}\label{sec:lattice}
A \emph{lattice} is a free $\Z$-module $L$ of finite rank
with a nondegenerate symmetric bilinear form
$\intfL{\phantom{i}, \phantom{i}}: L\times L\to \Z$,
which we call the \emph{intersection form}.
We let the group $\OG(L)$ of isometries of $L$ act on $L$ from the \emph{right},
and write the action of $g\in \OG(L)$ on $L$ by $x\mapsto x^g$.
A lattice $L$ is \emph{even} if $\intfL{v,v}$ is even for any vector $v\in L$.
A lattice $L$ is \emph{hyperbolic} if its rank $n$ is larger than $1$  and 
the real quadratic space $L\tensor\R$ is of signature $(1, n-1)$.
\par
Let $L$ be an even hyperbolic lattice.
The open subset $\shortset{x\in L\tensor\R}{\intfL{x,x}>0}$ of $L\tensor \R$
has two connected components,
each of which is called 
a \emph{positive cone}.
We choose a positive cone $\PPP_L$,
and denote by  $\OG^+(L)$ the stabilizer subgroup  of $\PPP_L$ in $\OG(L)$.
A vector $r\in L$ is called a \emph{$(-2)$-vector}
if $\intfL{r, r}=-2$.
Let $r$ be a $(-2)$-vector.
We put
$$
(r)\sperp :=\set{x\in \PPP_L}{\intfL{x, r}=0},
$$
and call it a \emph{$(-2)$-hyperplane}.
The reflection
$$
s_r: x\mapsto x+\intfL{x, r}\cdot r
$$ 
in $(r)\sperp$ is an element of $\OG^+(L)$.
We denote by  $\Weyl{L}$  the subgroup of $\OG^+(L)$
generated by all the reflections $s_r$ in $(-2)$-hyperplanes,
and call $\Weyl{L}$ the \emph{Weyl group} of $L$.
A \emph{standard fundamental domain of $\Weyl{L}$}
is the closure in $\PPP_L$ of a connected component of
$$
\PPP_L\;\setminus \;\bigcup_{r}\; (r)\sperp,
$$
where $r$ ranges through the set of $(-2)$-vectors.
Note that $\Weyl{L}$ acts on the set of standard fundamental domains
transitively.
\par
Suppose that a basis of an even hyperbolic lattice $L$
and the Gram matrix of the intersection form $\intfL{\phantom{i}, \phantom{i}}$ with respect to this basis
are given.
We have the following algorithms. See~\cite[Section 3]{MR3166075} for the details.
\begin{algorithm}\label{algo:affES}
 Let $v$ be a vector in  $\PPP_L \cap L$.
Then, for  an integer $a$ and  an even integer $d$, the finite set
$\shortset{x\in L }{ \intfL{x,v}=a, \intfL{x,x}=d}$ can be calculated.
In particular, the sets
$$
\RRR(v):=\set{r\in L}{\intfL{r, v}=0, \;\intfL{r, r}=-2}
$$ 
and 
$$
\FFF(v):=\set{f\in L}{\intfL{f, v}=1,\; \intfL{f, f}=0}
$$
can be calculated.
\end{algorithm}
\begin{algorithm}\label{algo:separating}
Let $u$ and $v$ be vectors in $\PPP_L\cap L$.
Then, for a negative even integer $d$, the finite set
$\shortset{x\in L }{ \intfL{x,u}>0, \intfL{x,v}<0, \intfL{x, x}=d}$ can be calculated.
In particular, the set 
$$
\SSS(u, v):=\set{r\in L }{ \intfL{r,u}>0, \;\intfL{r,v}<0, \;\intfL{r, r}=-2}
$$
can be calculated.
\end{algorithm}
We call an ordered nonempty set 
$$
\va:=[h_0, \rho_1, \dots, \rho_K]
$$
of vectors of $L$ an \emph{ample list of vectors}
if $h_0\in \PPP_L\cap L$ and,
for any $r\in \RRR(h_0)$, there exists a member $\rho_i$ of $\{  \rho_1, \dots, \rho_K\}$ such that $\intfL{r, \rho_i}\ne 0$.
\begin{example}
(1) If vectors $ \rho_1, \dots, \rho_K$ of $L$  span the linear space $L\tensor \Q$ over $\Q$, then 
$[h_0, \rho_1, \dots, \rho_K]$ is an ample list of vectors
for any vector $h_0\in \PPP_L\cap L$.
\par
(2)  If a vector $h_0\in \PPP_L\cap L$ satisfies  $\RRR(h_0)=\emptyset$, then the list $[h_0]$ is an ample list of vectors.
\par
(3) If $[h_0, \rho_1, \dots, \rho_K]$ is an ample list of vectors,
then $[h_0, \rho_1, \dots, \rho_K, \rho_{K+1}]$ is an ample list of vectors for any $\rho_{K+1}\in L$.

 \end{example}
 Let $\va=[h_0, \rho_1, \dots, \rho_K]$ be an ample list of vectors.
 %and let $\varepsilon$ be 
 %a parameter that moves in a region of  small positive real numbers.
 We define $D(\va)$ to be 
 the  unique standard  fundamental domain  of $\Weyl{L}$ such that 
$$
\vae:=h_0+\varepsilon \rho_1+\cdots+\varepsilon^K \rho_K
$$
is contained in the interior of $D(\va)$,
where $\varepsilon$  is a  sufficiently small positive real number.
For $x\in \PPP_L$, 
we write 
$$
\intfL{\va, x}>0
$$ 
if  the real vector
$$
(\;\intfL{h_0, x}, \, \intfL{\rho_1, x}, \,\dots, \, \intfL{\rho_K, x}\;) \in \R^{K+1}
$$
is nonzero and its leftmost nonzero entry
is positive;
that is, $\intfL{\vae, x}\in \R$  is positive for a sufficiently small positive real number $\varepsilon$.
For $x_1, x_2\in \PPP_L$, we write 
$$
\intfL{\va, x_1}>\intfL{\va, x_2}
$$ 
if $\intfL{\va, x_1-x_2}>0$.
We put
$$
\RRR^+ (\va):=\set{r\in \RRR(h_0)}{\intfL{\va, r}>0}.
$$
Note that $\RRR(h_0)$ is the disjoint union of $\RRR^+(\va)$ and $-\RRR^+(\va)$.
Then $D(\va)$ is the unique standard fundamental domain of $\Weyl{L}$ that contains $h_0$ and is contained in the region
$$
\set{x\in \PPP_L}{\intfL{x, r}\ge 0\;\;\textrm{for any vector}\;\; r\in \RRR^+ (\va)}.
$$
The following lemma is obvious.
\begin{lemma}\label{lem:inDa}
A vector $v\in \PPP_L\cap L$ is contained in $D(\va)$ if and only if
$\SSS(h_0, v)=\emptyset$ and $\intfL{v, r}\ge 0$ for any vector $r\in \RRR^+(\va)$.
\end{lemma}
Let $d$ be an even positive integer.
Suppose that a vector $v\in \PPP_L\cap L$ satisfies $\intfL{v, v}=d$.
From $v$, we can find a vector  $h_v$ in $D(\va)\cap L$ satisfying  $\intfL{h_v, h_v}=d$ by the following method.
First we calculate the union 
$$
\SSS(h_0, v) \cup \RRR\sprime =\{r_1, \dots, r_M\}, 
$$ 
where  %$\SSS(u, v\sprime)$ is defined in Algorithm~\ref{algo:separating}, and 
$$
\RRR\sprime:=\set{r\in \RRR^+(\va)}{\intfL{v, r}<0}.
$$
Note that we have  $\intfL{v, r_i}<0$  and $\intfL{\va, r_i}>0$
for each $r_i\in \SSS(h_0, v) \cup \RRR\sprime$.
Note also that, if a $(-2)$-vector $r$ satisfies $\intfL{v, r}<0$  and $\intfL{\va, r}>0$,
then $r$ belongs to $\SSS(h_0, v) \cup \RRR\sprime$.
We put
$$
\vt_i:=\frac{-1}{\intfL{v, r_i}} \left( \;\intfL{h_0, r_i}, \; \intfL{\rho_1, r_i}, \; \dots, \; \intfL{\rho_{K}, r_i}\; \right) \in \R^{K+1}.
$$
If $\vt_i=\vt_j$ holds  for some distinct indices $i$ and $j$,
then we choose a random vector $\rho_{K+1}\in L$ and replace $\va$ by a new ample list of vectors 
$$
[h_0, \rho_1, \dots, \rho_{K}, \rho_{K+1}].
$$
(Note that this replacement of $\va$ does not change $D(\va)$.)
Repeating this process,
we can assume that $\vt_1, \dots, \vt_M$ are  distinct.
We sort the vectors  $r_1, \dots, r_M$ of $\SSS(h_0, v) \cup \RRR\sprime$ in such a way that, 
if $i>j$,  then the leftmost  nonzero entry of $\vt_i-\vt_j$ is positive.
Consider the half-line $\ell$ in $\PPP_L$
given by
$$
\vae + tv \;\; (t\in \R_{\ge 0}),
$$
where $\varepsilon$ is a sufficiently small positive real number.
Then $\ell$ is not contained in any $(-2)$-hyperplane,
the $(-2)$-hyperplanes $(r_1)\sperp, \dots, (r_M)\sperp$
intersect $\ell$ at distinct points,
and any $(-2)$-hyperplane intersecting $\ell$ is one of $(r_1)\sperp, \dots, (r_M)\sperp$.
Moreover,  the values $t_i$ of the parameter $t$ of $\ell$ at which $\ell$ intersects $(r_i)\sperp$ satisfy 
$$
t_1>\dots >t_M>0, 
$$
because,  if $\vt_i=(t_{i, 0}, t_{i, 1}, \dots, t_{i, K})\in \R^{K+1}$, then we have
$$
t_i=t_{i, 0}+\varepsilon\, t_{i, 1}+ \dots + \varepsilon^K \,t_{i, K}.
$$
Therefore, if we denote by $s_i\in \Weyl{L}$ the reflection in $(r_i)\sperp$,
then the vector
\begin{equation}\label{eq:hvsss}
h_v:=v^{s_1\dots s_M}
\end{equation}
belongs to  $D(\va)\cap L$.
\section{Polarizations of degree $2$}\label{sec:poldeg2}
\begin{figure} 
\def\ha{40}
\def\hav{37}
\def\hd{25}
\def\hdv{22}
\def\he{10}
\def\hev{7}
\setlength{\unitlength}{1.2mm}
\vskip .5cm
\centerline{
{\small
\begin{picture}(100, 37)(-20, 7)
\put(0, \ha){$A\sb l$}
\put(10, \ha){\circle{1}}
\put(9.5, \hav){$a\sb 1$}
\put(10.5, \ha){\line(5, 0){5}}
\put(16, \ha){\circle{1}}
\put(15.5, \hav){$a\sb 2$}
\put(16.5, \ha){\line(5, 0){5}}
\put(22, \ha){\circle{1}}
\put(21.5, \hav){$a\sb 3$}
\put(22.5, \ha){\line(5, 0){5}}
\put(30, \ha){$\dots\dots\dots$}
\put(45, \ha){\line(5, 0){5}}
\put(50.5, \ha){\circle{1}}
\put(50, \hav){$a\sb {l}$}
\put(0, \hd){$D\sb m$}
\put(10, \hd){\circle{1}}
\put(9.5, \hdv){$d\sb 2$}
\put(10.5, \hd){\line(5, 0){5}}
\put(16, 31){\circle{1}}
\put(17.5, 30.5){$d\sb 1$}
\put(16, 25.5){\line(0,1){5}}
\put(16, \hd){\circle{1}}
\put(15.5, \hdv){$d\sb 3$}
\put(16.5, \hd){\line(5, 0){5}}
\put(22, \hd){\circle{1}}
\put(21.5, \hdv){$d\sb 4$}
\put(22.5, \hd){\line(5, 0){5}}
\put(30, \hd){$\dots\dots\dots$}
\put(45, \hd){\line(5, 0){5}}
\put(50.5, \hd){\circle{1}}
\put(50, \hdv){$d\sb {m}$}
\put(0, \he){$E\sb n$}
\put(10, \he){\circle{1}}
\put(9.5, \hev){$e\sb 2$}
\put(10.5, \he){\line(5, 0){5}}
\put(16, \he){\circle{1}}
\put(15.5, \hev){$e\sb 3$}
\put(22, 16){\circle{1}}
\put(23.5, 15.5){$e\sb 1$}
\put(22, 10.5){\line(0,1){5}}
\put(16.5, \he){\line(5, 0){5}}
\put(22, \he){\circle{1}}
\put(21.5, \hev){$e\sb 4$}
\put(22.5, \he){\line(5, 0){5}}
\put(30, \he){$\dots\dots\dots$}
\put(45, \he){\line(5, 0){5}}
\put(50.5, \he){\circle{1}}
\put(50, \hev){$e\sb {n}$}
\end{picture}
}
}
\vskip 10pt
\caption{ Indecomposable $ADE$-configurations}\label{fig:ADE}
\end{figure}
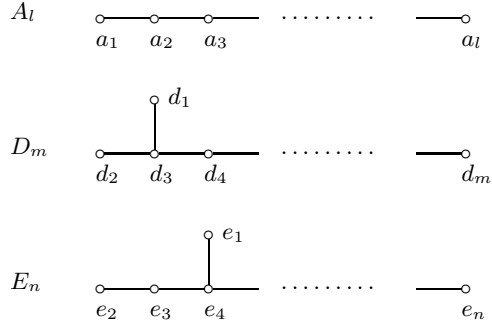

Let $X$ be a $K3$ surface defined over an algebraically closed field $k$ of  characteristic not equal to $2$,
and let $S_X$ denote the N\'eron--Severi lattice of $X$ 
with the intersection form $\intfS{\phantom{i}, \phantom{i}}$.
Suppose that  $\rank S_X$ is larger than $1$.
Then $S_X$ is an even hyperbolic lattice.
We let the automorphism group $\Aut(X)$ act on $X$ from the left and act on $S_X$ from the right by the pull-back.
Let $\PPP(X)$ denote the positive cone of $S_X$ that contains  an ample class.
We put
$$
N(X):=\set{x\in \PPP(X)}{\intfS{x, [C]}\ge 0\;\;\textrm{for any curve}\;\; C\subset X},
$$
where $[C]\in S_X$ is the class of a curve $C$ on $X$.
%and let $N(X)\spcirc$ denote the interior of $N(X)$.
It is well known that $N(X)$ is a standard fundamental domain of the Weyl group $\Weyl{S_X}$.
A vector $h\in S_X$ with $\intfS{h, h}=2$ is called a \emph{polarization of degree $2$}
if the complete linear system $|\LLL_h|$ of a line bundle $\LLL_h\to X$ whose class  is $h$
is fixed-component free.
By~\cite{MR1260944}, 
we have the following criterion.
\begin{proposition}\label{prop:Nikulindeg2pol}
A vector $h\in S_X$ with $\intfS{h, h}=2$ is  a polarization of degree $2$
if and only if $h\in N(X)$ and $\FFF(h)=\emptyset$. %, where $\FFF(h)$ is defined in Algorithm~\ref{algo:separating}.
\end{proposition}
Suppose that $h\in S_X$ is a polarization of degree $2$.
Then, by~\cite{MR0364263},  the complete linear system $|\LLL_h|$ is base-point free, and hence defines a generically finite morphism 
$\Phi_h: X\to \P^2$ of degree $2$.
Let 
$$
X\;\maprightsp{\psi_h}\; Y_h \;\maprightsp{\pi_h} \;\P^2
$$
be the Stein factorization of $\Phi_h$,
and let $B_h\subset \P^2$ be the branch curve of the double covering $\pi_h$.
Then $\psi_h: X\to Y_h$ is a contraction of smooth rational curves, and $B_h$ is a curve of degree $6$ with only simple singularities.
For each singular point $P$ of $B_h$,
the  curves contracted to $P$ by $\Phi_h$
form an indecomposable $ADE$-configuration of smooth rational curves.
We put
$$
\EEE_P(h)\;:=\; \set{[C]}{\textrm{$C$ is a smooth rational curve on $X$ contracted to $P$ by $\Phi_h$}},
$$
%Then the vectors of $\EEE_P(h)$ form an indecomposable 
%root system of type $A_l$ ($l\ge 1$), $D_m$ ($m\ge 4$),  or $E_n$ ($n=6,7,8$).
%$ADE$-configuration of smooth rational curves on $X$.
and label the elements of $\EEE_P(h)$ in such a  way that their dual graph is indicated  in Figure~\ref{fig:ADE}.
\par
We denote by $\tau(h)\in\Aut(X)$ the involution of $X$
induced by  the Galois transformation of 
the double covering $\pi_h$,
and call it a \emph{double plane involution}.
Suppose that a basis of $S_X$ and the Gram matrix of $\intfS{\phantom{i}, \phantom{i}}$
with respect to this basis are given.
Suppose also that 
we have an ample list of vectors $\va$ such that
$$
D(\va)=N(X)
$$
holds.
Then we can calculate the matrix representation $M(h)$ 
of the action 
of $\tau(h)$ on $S_X$   by the following method.
It is well known that there exists a successive blowing up $\beta_h: F_h\to \P^2$ of $\P^2$
at  (possibly infinitely near) points of the singular locus of $B_h$
such that
$\Phi_h$ factors as 
$$
X\maprightsp{q_h} F_h\maprightsp{\beta_h} \P^2,
$$
where $q_h$ is the quotient morphism by $\tau(h)$.
Let $S_F$ denote the N\'eron--Severi lattice of the smooth  rational surface  $F_h$.
Then the pull-back $q_h^*$ by $q_h$ identifies $S_F\tensor \Q$ with the eigenspace of $\tau(h)$ in $S_X\tensor\Q$
with eigenvalue $1$,
and hence $\tau(h)$ acts on the orthogonal complement of $q_h^*S_F\tensor\Q$ in  $S_X\tensor\Q$ as
the scalar multiplication by $-1$.
On the other hand, the subspace $q_h^*S_F\tensor\Q$ 
is generated by $h$  and %, which is the pull-back of the class of a line on $\P^2$,  and
the vectors of the form $r+r^{\tau(h)}$,
where $r\in \EEE_P(h)$ and $P\in \Sing (B_h)$.
The action of $\tau(h)$  on $\EEE_P(h)$ is as follows:
\begin{itemize}
\item If $P$ is of type $A_{l}$, then $a_i^{\tau(h)}=a_{l+1-i}$ for $i=1, \dots, l$.
\item If $P$ is of type $D_{2k}$, then $\tau(h)$ acts on $\EEE_P(h)$ as the identity.
\item If $P$ is of type $D_{2k+1}$, then $d_1^{\tau(h)}=d_{2}$, $d_2^{\tau(h)}=d_{1}$,
and $d_i^{\tau(h)}=d_{i}$  for $i=3, \dots, 2k+1$.
\item If $P$ is of type $E_6$, then $e_1^{\tau(h)}=e_{1}$, and $e_i^{\tau(h)}=e_{8-i}$ for $i=2, \dots, 6$.
\item If $P$ is of type $E_7$ or $E_8$,  then $\tau(h)$ acts on $\EEE_P(h)$ as the identity.
\end{itemize}
Hence,  in order to calculate 
the matrix representation $M(h)$ of $\tau(h)$ on $S_X$, 
it is enough to calculate the sets $ \EEE_P(h)$.
\par
We put
$$
\EEE(h):=\bigcup_{P\in \Sing(B_h)} \EEE_P(h).
$$
First we calculate the finite set 
$$
\RRR^+(h):=\set{r\in \RRR(h)}{\intfS {\va, r}>0}.
$$
%where $ \RRR(h)$ is defined in Algorithm~\ref{algo:affES}.
Note that, since  $D(\va)$ is equal to $N(X)$ and any $r\in \EEE(h)$ is the class of a curve, 
we have $\intfS{\va, r}>0$ for any vector $r\in \EEE(h)$.
Moreover, any vector $r\sprime\in \RRR^+(h)$ is the class of an effective divisor,
each irreducible component of which is a smooth rational curve contracted by $\Phi_h$.
Therefore, we have $\EEE(h)\subset \RRR^+(h)$.
Moreover, 
a vector $r\sprime\in \RRR^+(h)$ is 
a linear combination with nonnegative integer coefficients of vectors in $\EEE(h)$.
Consequently,  a vector 
$r\sprime\in \RRR^+(h)$ does  \emph{not} belong to $\EEE(h)$
if and only if
$r\sprime$ can be written as a linear combination with nonnegative integer coefficients of vectors $r\spprime$ in $ \RRR^+(h)$
satisfying  $\intfS{\va, r\spprime}< \intfS{\va, r\sprime}$.
Thus, starting from the vector $r_0$ of $\RRR^+(h)$ with the smallest $\intfS{\va, r_0}$,
we can successively detect the elements of $\EEE(h)$ in $\RRR^+(h)$.
We connect two distinct elements $r, r\sprime$ of $\EEE(h)$ by an edge if and only if
$\intfS{r, r\sprime}=1$.
Then 
the vertices of each  connected component of $\EEE(h)$ form the set  $\EEE_P(h)$.
\begin{remark}
This method of calculating the action of $\tau(h)$ on $S_X$ was also used in finding  a finite set of generators
of  $\Aut(X)$ by Borcherds method 
in~\cite{MR3190354} and~\cite{1412.6904},
and in the study of projective models of the supersingular $K3$ surface $X(5)$
in characteristic $5$ with Artin invariant $1$ in~\cite{MR3166075}.
\end{remark}
\section{N\'eron--Severi lattices of supersingular $K3$ surfaces}\label{sec:RS}
%
%Let $X$ be a supersingular $K3$ surface with Artin invariant $\sigma$ in \emph{odd} characteristic $p$.
Rudakov and Shafarevich~\cite{MR633161} 
proved the following theorems. 
For the proof of Theorem~\ref{thm:RS1}, 
see also~\cite[Chapter 15]{MR1662447}.
\begin{theorem}\label{thm:RS1}
Let $p$ be an odd prime, and let $\sigma$ be  a positive integer less than or equal to $10$.
Then there exists a lattice $\Lpsm$, unique up to  isomorphism,
with the following properties.
{\rm (i)} $\Lpsm$ is an even hyperbolic lattice of rank $22$.
{\rm (ii)} The discriminant group $(\Lpsm)\dual/\Lpsm$ of $\Lpsm$ is isomorphic to $(\Z/p\Z)^{2\sigma}$.
\end{theorem}
\begin{theorem}\label{thm:RS2}
Let $X$ be a supersingular $K3$ surface in odd characteristic $p$ with Artin invariant $\sigma$.
Then its N\'eron--Severi lattice $S_X$ is isomorphic to $\Lambda_{p, \sigma}^{-}$.
\end{theorem}
An explicit method of constructing  $\Lpsm$ is also given in~\cite{MR633161} 
(see also~\cite{MR2036331}).
We use the following construction,
which is slightly different from the one given in~\cite{MR633161}.
The ingredients of the construction are the following  lattices.
\par
(i) Let $U$ and $\Up$ be the even hyperbolic lattices of rank $2$ with the Gram matrices
\begin{equation}\label{eq:GramUUp}
\left[\begin{array}{cc} 0 & 1 \\ 1 & 0\end{array}\right]
\quand 
\left[\begin{array}{cc} 0 & p \\ p & 0\end{array}\right],
\end{equation}
respectively.
\par
(ii) Let $q$ be a prime satisfying
$$
q\equiv 3 \bmod 8 \quand  \left(\frac{-q}{p}\right)=-1, 
$$
and let $\gamma$ be an integer satisfying $\gamma^2+p\equiv 0 \bmod q$.
Let $\Hp$ be the even \emph{negative} definite lattice of rank $4$
with the Gram matrix
$$
\renewcommand{\arraystretch}{1.4}
(-1)\left[\begin{array}{cccc} 
2 & 1 & 0 & 0 \\
1 & (q+1)/2 & 0 & \gamma\\
0 & 0  & p(q+1)/2 & p \\
0 & \gamma & p & 2(p+\gamma^2)/q
\end{array}\right].
$$
Then the discriminant group of $H^{(-p)}$ is isomorphic to $(\Z/p\Z)^2$.
See~\cite{MR0568309} and~\cite{MR2036331}.
\par
(iii) 
Let $E_8$ denote the root lattice of type $E_8$,
which is an even unimodular positive definite lattice of rank $8$.
Then $E_8$ has a \emph{standard basis }$e_1, \dots, e_8$,
whose dual graph is given in Figure~\ref{fig:ADE}.
Let $\mE$ be the lattice obtained from $E_8$ by multiplying the intersection form by $-1$,
and let $\pmE$ be the lattice obtained from $\mE$ by multiplying the intersection form by $p$.
Then the discriminant group of $\pmE$ is isomorphic to $(\Z/p\Z)^8$.
\par
Then $\Lpsm$ is isomorphic to the following lattices:
\begin{equation*}\label{eq:UHEE}
\renewcommand{\arraystretch}{1.4}
\begin{array}{ll}
U\oplus \Hp \oplus \mE \oplus \mE & \textrm{if $\sigma=1$, }\\
\Up\oplus \Hp \oplus  \mE \oplus \mE & \textrm{if $\sigma=2$, }\\
U\oplus \Hp \oplus  \Hp \oplus \Hp  \oplus \mE & \textrm{if $\sigma=3$, }\\
\Up\oplus \Hp \oplus  \Hp \oplus \Hp   \oplus \mE & \textrm{if $\sigma=4$, }\\
U\oplus \Hp  \oplus \mE \oplus \pmE & \textrm{if $\sigma=5$, }\\
\Up\oplus \Hp \oplus \mE \oplus \pmE & \textrm{if $\sigma=6$, }\\
U\oplus  \Hp \oplus  \Hp \oplus \Hp  \oplus \pmE & \textrm{if $\sigma=7$, }\\
\Up\oplus \Hp \oplus  \Hp \oplus \Hp   \oplus \pmE & \textrm{if $\sigma=8$, }\\
U\oplus \Hp \oplus  \pmE \oplus \pmE & \textrm{if $\sigma=9$, }\\
\Up\oplus \Hp \oplus  \pmE \oplus \pmE & \textrm{if $\sigma=10$.}\\
\end{array}
%\label{eq:4cases}
\end{equation*}
Let $\intfLL{\phantom{i}, \phantom{i}}$ denote the intersection form of $\Lpsm$.
Note that $\Lpsm$ has the form of the orthogonal direct sum 
$$
U\sprime\oplus N,
$$
where $U\sprime$ is $U$ or $\Up$ according to the parity of $\sigma$,
and $N$ is an even negative definite lattice 
with the intersection form $\intfN{\phantom{i}, \phantom{i}}$.
We put
$$
p\sprime:=\begin{cases}
1 & \textrm{if $U\sprime$ is $U$}, \\
p & \textrm{if $U\sprime$ is $\Up$}.
\end{cases}
$$
We choose a vector $n\in N$ randomly.
If $2-\intfN{n, n}$ is divisible by $2\,p\sprime$,
then we can find a vector $u\in U\sprime$ such that $v:=u+n\in \Lpsm$ satisfies 
$\intfLL{v, v}=2$.
By this method, we can generate many vectors of $\Lpsm$ with square-norm $2$.
\section{Generating double plane involutions}
We fix an odd prime $p$ and a positive integer $\sigma$ less than or equal to $10$.
Let $X$ be a supersingular $K3$ surface in characteristic $p$
with Artin invariant $\sigma$.
We make a set $\MMM$ of matrix representations 
on $S_X$ of double plane involutions $\tau(h)\in \Aut(X)$ associated with
 polarizations $h\in S_X$ of degree $2$.
 \begin{itemize}
 \item[(0)] We set $\MMM=\{\}$.
 \item[(1)] We construct a Gram matrix of the lattice $\Lpsm$ by the result in Section~\ref{sec:RS}.
 \item[(2)] We find a vector $h_0\in \Lpsm$ such that $\intfLL{h_0, h_0}>0$.
Let  $\PPP_{\Lambda}$ be the positive cone of $\Lpsm$
 containing $h_0$.
\item[(3)]
 We calculate $\RRR(h_0)$,
 and choose an ample list of vectors
 $$
 \va:=[h_0, \rho_1, \dots, \rho_K].
 $$
\item[(4)]
 By Theorem~\ref{thm:RS2},
there exists an isomorphism $\iota: \Lpsm\isom S_X$ of lattices.
 Multiplying $\iota$ by $-1$ if necessary,
 we can assume that $\iota$ maps $\PPP_{\Lambda}$ to $\PPP(X)$.
 Composing $\iota$ with an element of $\Weyl{S_X}$  if necessary,
 we can further assume that $\iota$ maps $D(\va)$ to $N(X)$.
 From now on, we identify $\Lpsm$ with $S_X$,  and $D(\va)$ with $N(X)$ by 
 the isometry $\iota$.
 \item[(5)]
 We make a finite  set $\VVV$ of vectors $v\in \Lpsm$ with $\intfLL{v, v}=2$
 by the method described in Section~\ref{sec:RS}.
\item[(6)] For each $v\in \VVV$, we execute  the following calculations.
\begin{itemize}
\item[(6-1)] If $\intfLL{v, h_0}<0$, then we replace $v$ with $-v$,
 so that we can assume that $v\in \PPP_{\Lambda}$.
\item[(6-2)]  We calculate $\FFF(v)$.
 If $\FFF(v)\ne\emptyset$, we proceed  to the next element of $\VVV$.
 If $\FFF(v)=\emptyset$, we go to Step (6-3).
\item[(6-3)]
 From $v$, we construct the vector $h_v\in \Lpsm$  with $\intfLL{h_v, h_v}=2$ that belongs to $D(\va)$
 by the method described in Section~\ref{sec:lattice}.
 Since $h_v$ and $v$ are related by~\eqref{eq:hvsss},
 we have  $\FFF(h_v)=\emptyset$.
By the identification of $D(\va)$ with $N(X)$,
 we see that $h_v$ is nef.
 Therefore, by Proposition~\ref{prop:Nikulindeg2pol},
 we see that $h_v$ is a polarization of degree $2$.
\item[(6-4)]
 We then calculate the matrix representation $M(h_v)$ of 
 the double plane involution $\tau(h_v)\in \Aut(X)$ by the method described in Section~\ref{sec:poldeg2},
 and append $M(h_v)$ to $\MMM$.
 \end{itemize}
  \end{itemize}
Once we make a sufficiently large set
$$
\MMM=\{M(h_1), \dots, M(h_N)\}
$$
of $22\times 22$ matrices representing the action of double plane involutions of $X$ on $S_X$,
we make a product
$$
M:=M(h_{i_1}) \cdots M(h_{i_\nu})
$$
of randomly chosen elements of $\MMM$,
and calculate its characteristic polynomial $\phi_{M}(t)$.
By Theorem~\ref{thm:cycsSalem},
if $\phi_M(t)$ is irreducible in $\Z[t]$ and not equal to the cyclotomic polynomial $(t^{23}-1)/(t-1)$,
then $\phi_M(t)$ is a Salem polynomial.
\par
By this method,
we confirm that, if $p$ is an odd prime $\le \thep$,
then $\Aut(X)$ contains an automorphism of irreducible Salem type
that is a product of at most $22$ double plane involutions.
\begin{remark}\label{rem:infample}
Let $\ve_1, \dots, \ve_{22}$ be a basis of $\Lpsm$,
and let  $\ve_1\dual, \dots, \ve_{22}\dual$ be the dual basis.
Note that $p\ve_i\dual \in \Lpsm$ holds for $i=1, \dots, 22$.
Hence, in Step (3), 
 we can choose 
$[h_0, p\ve_1\dual, \dots, p\ve_{22}\dual]$
as an ample list of vectors.
\end{remark}
\section{Supersingular $K3$ surfaces with Artin invariant $10$}\label{sec:sigma10}
We consider a supersingular $K3$ surface $X$ in characteristic $p\ge 11$ with Artin invariant $10$.
We have  
$$
\Lambda_{p, 10}^{-}=\Up\oplus \Hp\oplus \pmE \oplus \pmE.
$$
Let $u_1, u_2$ be the basis of $\Up$ with the Gram matrix ~\eqref{eq:GramUUp},
and let $e_1, \dots, e_8$ (resp.~$e_1\sprime, \dots, e_8\sprime$) be the standard basis of the first $\pmE$ (resp.  the second $\pmE$).
 In particular, each $e_{\nu}$ or $e_{\nu}\sprime$ is of square-norm $-2p$.
 For $v\in \Hp$ and $a\in \Z$,
 we denote by
 $$
 (a, 1, v)\; \in\;  \Up\oplus\Hp
 $$
 the vector $au_1+u_2+v$.
 Then the square-norm of $(a, 1, v)$ is $2pa+\intfH{v, v}$,
 where  $\intfH{\phantom{a}, \phantom{a}}$ is the intersection form of $\Hp$.
 Note that, if $(a, 1, v)\in \Up\oplus\Hp$ is of square-norm $2$, then  the vectors $(a+1, 1, v)+e_{\nu}$ and $(a+1, 1, v)+e\sprime_{\nu}$
 of $\Lambda_{p, 10}^{-}$ are also of square-norm $2$
 for $\nu=1, \dots, 8$.
\par
\begin{figure}
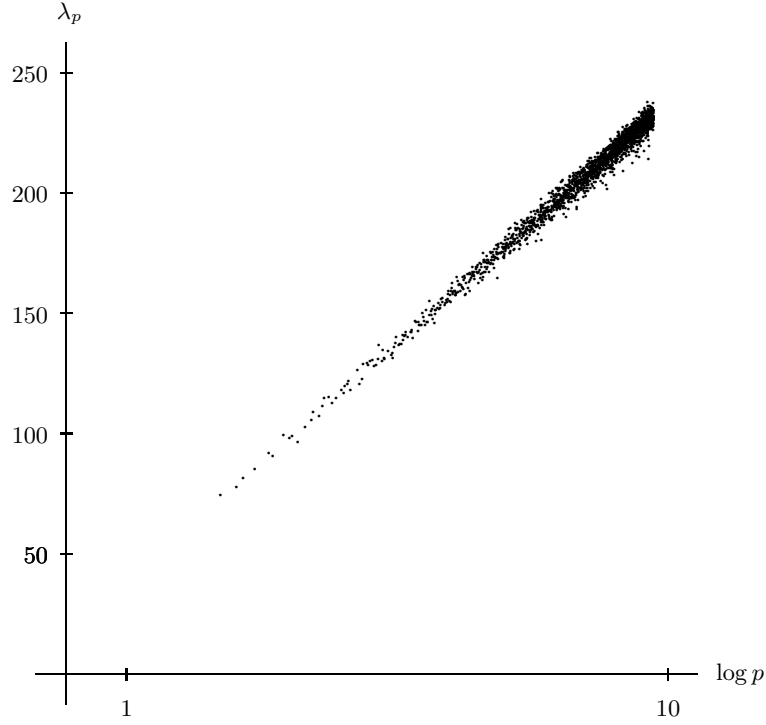

{\small
\setlength{\unitlength}{.8cm}

}
\caption{Growth of the entropy}\label{fig:entropy}
\end{figure}

For $p$ with $11\le p\le \thepten$, we have found  six vectors $v_k \in \Hp$ and six positive integers $a_k\in \Z$
with the following properties (i)--(v).
\begin{itemize}
\item[(i)] The vector $h_k:=(a_k, 1, v_k)$ is of square-norm $2$ for $k=1, \dots, 6$.
\end{itemize}
We put
$$
h_{6+\nu} :=(a_k+1, 1, v_k)+e_{\nu}, \quad
h_{14+\nu} :=(a_k+1, 1, v_k)+e\sprime_{\nu},
$$
for $\nu=1, \dots, 8$.
Then $h_7, \dots, h_{22}$ are also of square-norm $2$.
\begin{itemize}
\item[(ii)] $\intfLL{h_1, h_i}>0$ for $i=2, \dots, 22$.
\item[(iii)] $\SSS(h_1, h_i)=\emptyset$ for $i=2, \dots, 22$.
\item[(iv)] $\RRR(h_i)=\emptyset$ and $\FFF(h_i)=\emptyset$ for $i=1, \dots, 22$.
\end{itemize}
Since $R(h_1)=\emptyset$,
there exists a unique standard fundamental domain $D([h_1])$ of the Weyl group  $W(\Lpsm)$ 
that contains $h_1$ in its interior.
Since  $\SSS(h_1, h_i)=\emptyset$ for $i=2, \dots, 22$,
we see that $h_1, \dots, h_{22}$ are also contained in $D([h_1])$.
Hence, under a suitable isometry $\Lambda_{p, 10}^{-}\isom S_X$,
we can assume that each $h_i$ is a nef vector in $S_X$.
Since $\FFF(h_i)=\emptyset$ for $i=1, \dots, 22$, 
we see that each $h_i$ is 
a polarization of degree $2$ on $X$.
Moreover, since $\RRR(h_i)=\emptyset$,
the branch curve $B_{h_i}\subset \P^2$ 
of the double plane involution $\tau(h_i)$ is smooth.
Hence  $\tau(h_i)$ acts on $h_i$ trivially,
and on the orthogonal compliment of $h_i$ as the multiplication by $-1$.
\begin{itemize}
\item[(v)] The product $g:=\tau(h_1)\cdots \tau (h_{22})$ is of irreducible Salem type.
\end{itemize}
This observation and a computer-aided calculation give the proof of Theorem~\ref{thm:main2}.
\begin{example}
Consider the case $p=\thepten$.
Then $\Hp$ has a Gram matrix
$$
\left[ \begin {array}{cccc} -2&-1&0&0\\ \noalign{\medskip}-1&-30&0&-4
\\ \noalign{\medskip}0&0&-521670&-17389\\ \noalign{\medskip}0&-4&-
17389&-590\end {array} \right]
$$
under a certain basis $\eta_1, \dots, \eta_4$ of $\Hp$.
The vectors
%   [1, 1, 15, 31, 0, -3]  [1, 1, 9, 18, -1, 25][1, 1, 51, 4, 0, -7]  [1, 1, 30, 29, 0, 3]   [1, 1, 55, -4, 0, 7] [2, 1, 19, 23, -2, 56]
%
$$
\begin{array}{rccrrrrrrrc}
h_1 &=& [&1, &1, &15, &31, &0, &-3&],\\
h_2 &=&[&1, &1, &9, &18, &-1, &25&],\\
h_3 &=&[&1, &1, &51, &4, &0, &-7&],\\
h_4 &=&[&1, &1, &30, &29, &0, &3&],\\
h_5 &=&[&1, &1, &55, &-4, &0, &7&],\\
h_6 &=&[&2, &1, &19, &23, &-2, &56&],
\end{array}
$$
of $\Up\oplus \Hp$ written with respect to the basis $u_1, u_2, \eta_1, \dots, \eta_4$
satisfies the properties (i)--(v).
The characteristic polynomial on $S_X$  of the automorphism $g$
obtained from these six vectors has a real root 
${4.2539 \dots} \times 10^{100}$.
%read "compdatas/sigma10IS17.txt":
%SalemNumb[ithprime(2000), 10];
%4.253918000*10^100
\end{example}
\begin{remark}
Let $g_p$ be the automorphism of a supersingular $K3$ surface $X$ with Artin invariant $10$
in characteristic $p$
obtained by the method described in this section,
let $\rho_p$ be the real root $>1$ of the characteristic polynomial of $g_p$ on $S_X$, 
and let $\lambda_p:=\log \rho_p$ be the \emph{entropy} of $g_p$.
Then, for $11\le p\le \thepten$,
we have
$$
\lambda_p \sim 19.1+ 21.8\, \log p.
$$
See Figure~\ref{fig:entropy}.
\end{remark}
\section{An example with Artin invariant $1$}\label{sec:example}
We denote by  $X(p)$ a supersingular $K3$ surface in characteristic $p$ with Artin invariant $1$,
which is unique up to isomorphism by the result of Ogus~\cite{MR563467, MR717616}. 
The existence of an automorphism $g\in \Aut(X(p))$ of irreducible Salem type 
was established by Blanc and Cantat~\cite{1307.0361} for $p=2$,  by Esnault and Oguiso~\cite{1406.2761} for $p=3$, and 
by Esnault, Oguiso, and Yu ~\cite{1411.0769} for $p=11$ or $p \ge 17$.
On the other hand, in~\cite{1502.06923}, Sch\"utt showed that, 
if  $p$ is odd and  satisfies $p \equiv 2 \bmod 3$, then 
there exists  a non-liftable automorphism  of  $X(p)$
whose characteristic polynomial  on  $S_{X(p)}$
is divisible by  a Salem polynomial of degree $20$.
\par
We consider the supersingular $K3$ surface $X(7)$,
which has not yet been  treated  by the previous works.
The lattice  $\Lambda_{7, 1}^{-}=U\oplus H^{(-7)} \oplus \mE\oplus \mE$ has a basis 
$\ve_1, \dots, \ve_{22}$ such that 
$\ve_1$ and $\ve_2$ form a basis of $U$ with the Gram matrix~\eqref{eq:GramUUp},
$\ve_3, \dots, \ve_6$  form a basis of $H^{(-7)}$ with the Gram matrix
$$
\left[ \begin {array}{cccc} -2&-1&0&0\\ \noalign{\medskip}-1&-6&0&-2
\\ \noalign{\medskip}0&0&-42&-7\\ \noalign{\medskip}0&-2&-7&-2
\end {array} \right] , 
$$
and $\ve_7, \dots, \ve_{14}$  (resp.,  $\ve_{15}, \dots, \ve_{22}$) form the standard basis of the first $\mE$ (resp., the second $\mE$).
We put
$$
h_0 := [1,1,0,0,0,0,0,0,0,0,0,0,0,0,0,0,0,0,0,0,0,0] \in \Lambda_{7, 1}^{-},
$$
which is of square-norm $2$.
The set $\RRR(h_0)$ consists of $486$ vectors.
The list
$$
\va:=[h_0, 7\ve_1\dual , \dots, 7\ve_{22}\dual]
$$
is an ample list of vectors.
We identify $\Lambda_{7, 1}^{-}$ with $S_{X(7)}$ by an isometry $\Lambda_{7, 1}^{-}\isom S_{X(7)}$
that maps $D(\va)$ to $N(X(7))$.
(Since $\FFF(h_0)\ne \emptyset$,
the vector $h_0$ is \emph{not} a polarization of degree $2$.)
\par
We consider the three vectors
%
%[5, 5, -2, 3, 2, -11, -12, -8, -16, -24, -20, -15, -10, -5, -8,  -5, -10, -15, -12, -9, -6, -3]
%[5, 5, -1, 0, 0, -2, -13, -9, -17, -25, -20, -15, -10, -5, -11,  -7, -14, -21, -17, -13, -9, -5]
%[3, 6, -2, 2, 2, -9, -5, -4, -7, -10, -8, -6, -4, -2, 0, 0, 0, 0,  0, 0, 0, 0]
%
\begin{eqnarray*}
h_1 &:=&[5, 5, -2, 3, 2, -11, -12, -8, -16, -24, -20, -15, -10,\\
&& \qquad  -5, -8,  -5, -10, -15, -12, -9, -6, -3], \\
h_2 &:=&[5, 5, -1, 0, 0, -2, -13, -9, -17, -25, -20, -15, -10, \\
&& \qquad  -5, -11,  -7, -14, -21, -17, -13, -9, -5],\\
h_3 &:=&[3, 6, -2, 2, 2, -9, -5, -4, -7, -10, -8, -6, -4, -2, 0, 0, 0, 0,  0, 0, 0, 0],
\end{eqnarray*}
of square-norm $2$.
 By means of Lemma~\ref{lem:inDa},
 we can confirm that $h_1, h_2, h_3$ are located in $D(\va)=N(X(7))$.
 Moreover we have $\FFF(h_1)=\FFF(h_2)=\FFF(h_3)=\emptyset$.
 Hence these $h_i$ are polarizations of degree $2$,
 and induce double plane involutions $\tau(h_i)$.
 The type of the singularities of the branch curve $B_{h_i}$ is
 %
 %
 %                   A[7] + A[4] + A[5]                  A[9] + 2 A[1] + A[7]                   E[8] + D[7] + A[2]
 %
 $$
A_4+A_5+A_7, \quad
2A_1+A_7+A_9, \quad
A_2+D_7+E_8, 
 $$
 respectively.
\begin{figure}
{\tiny 
\setlength{\arraycolsep}{1.8pt} 
$$
\left[ \begin {array}{cccccccccccccccccccccc} 24&24&-10&15&10&-55&-57
&-38&-76&-114&-95&-71&-48&-24&-40&-25&-50&-75&-60&-45&-30&-15
\\ \noalign{\medskip}24&24&-10&15&10&-55&-57&-38&-76&-114&-95&-72&-48&
-24&-40&-25&-50&-75&-60&-45&-30&-15\\ \noalign{\medskip}5&5&-3&3&2&-11
&-12&-8&-16&-24&-20&-15&-10&-5&-8&-5&-10&-15&-12&-9&-6&-3
\\ \noalign{\medskip}30&30&-12&17&12&-66&-72&-48&-96&-144&-120&-90&-60
&-30&-48&-30&-60&-90&-72&-54&-36&-18\\ \noalign{\medskip}-35&-35&14&-
21&-15&77&84&56&112&168&140&105&70&35&56&35&70&105&84&63&42&21
\\ \noalign{\medskip}10&10&-4&6&4&-23&-24&-16&-32&-48&-40&-30&-20&-10&
-16&-10&-20&-30&-24&-18&-12&-6\\ \noalign{\medskip}0&0&0&0&0&0&0&1&0&0
&0&0&0&0&0&0&0&0&0&0&0&0\\ \noalign{\medskip}0&0&0&0&0&0&1&0&0&0&0&0&0
&0&0&0&0&0&0&0&0&0\\ \noalign{\medskip}0&0&0&0&0&0&0&0&0&1&0&0&0&0&0&0
&0&0&0&0&0&0\\ \noalign{\medskip}0&0&0&0&0&0&0&0&1&0&0&0&0&0&0&0&0&0&0
&0&0&0\\ \noalign{\medskip}4&5&-2&3&2&-11&-10&-7&-14&-20&-16&-12&-8&-4
&-8&-5&-10&-15&-12&-9&-6&-3\\ \noalign{\medskip}1&-1&0&0&0&0&0&0&0&0&0
&0&0&0&0&0&0&0&0&0&0&0\\ \noalign{\medskip}0&1&0&0&0&0&-3&-2&-4&-6&-5&
-4&-3&-2&0&0&0&0&0&0&0&0\\ \noalign{\medskip}0&0&0&0&0&0&0&0&0&0&0&0&0
&1&0&0&0&0&0&0&0&0\\ \noalign{\medskip}5&5&-2&3&2&-11&-12&-8&-16&-24&-
20&-15&-10&-5&-9&-6&-12&-18&-15&-12&-8&-4\\ \noalign{\medskip}0&0&0&0&0
&0&0&0&0&0&0&0&0&0&0&0&0&0&0&0&0&1\\ \noalign{\medskip}0&0&0&0&0&0&0&0
&0&0&0&0&0&0&0&0&0&0&0&0&1&0\\ \noalign{\medskip}0&0&0&0&0&0&0&0&0&0&0
&0&0&0&0&0&0&0&0&1&0&0\\ \noalign{\medskip}0&0&0&0&0&0&0&0&0&0&0&0&0&0
&0&0&0&0&1&0&0&0\\ \noalign{\medskip}0&0&0&0&0&0&0&0&0&0&0&0&0&0&0&0&0
&1&0&0&0&0\\ \noalign{\medskip}0&0&0&0&0&0&0&0&0&0&0&0&0&0&0&0&1&0&0&0
&0&0\\ \noalign{\medskip}0&0&0&0&0&0&0&0&0&0&0&0&0&0&0&1&0&0&0&0&0&0
\end {array} \right] 
$$
}
\caption{$M(h_1)$}\label{fig:Mh1}
\end{figure}
\begin{figure}
{\tiny
\setlength{\arraycolsep}{2.4pt} 
$$
\left[ \begin {array}{cccccccccccccccccccccc} 6&6&0&0&0&-3&-15&-10&-
20&-29&-24&-18&-12&-6&-12&-8&-16&-24&-20&-16&-12&-6
\\ \noalign{\medskip}6&6&0&0&0&-3&-15&-10&-20&-30&-24&-18&-12&-6&-12&-
8&-16&-24&-20&-16&-12&-6\\ \noalign{\medskip}0&0&1&0&0&0&0&0&0&0&0&0&0
&0&0&0&0&0&0&0&0&0\\ \noalign{\medskip}6&6&1&-1&0&-2&-16&-10&-20&-30&-
24&-18&-12&-6&-12&-8&-16&-24&-20&-16&-12&-6\\ \noalign{\medskip}21&21&0
&0&-1&-7&-56&-35&-70&-105&-84&-63&-42&-21&-42&-28&-56&-84&-70&-56&-42&
-21\\ \noalign{\medskip}6&6&0&0&0&-3&-16&-10&-20&-30&-24&-18&-12&-6&-
12&-8&-16&-24&-20&-16&-12&-6\\ \noalign{\medskip}0&1&0&0&0&-1&0&0&0&0&0
&0&0&0&0&0&0&0&0&0&0&0\\ \noalign{\medskip}0&0&0&0&0&0&0&0&0&0&0&0&0&0
&0&0&0&0&0&0&0&1\\ \noalign{\medskip}0&1&0&0&0&0&0&0&0&0&0&0&0&0&-3&-2
&-4&-6&-5&-4&-3&-2\\ \noalign{\medskip}1&-1&0&0&0&0&0&0&0&0&0&0&0&0&0&0
&0&0&0&0&0&0\\ \noalign{\medskip}0&1&0&0&0&0&-3&-2&-4&-6&-5&-4&-3&-2&0
&0&0&0&0&0&0&0\\ \noalign{\medskip}0&0&0&0&0&0&0&0&0&0&0&0&0&1&0&0&0&0
&0&0&0&0\\ \noalign{\medskip}0&0&0&0&0&0&0&0&0&0&0&0&1&0&0&0&0&0&0&0&0
&0\\ \noalign{\medskip}0&0&0&0&0&0&0&0&0&0&0&1&0&0&0&0&0&0&0&0&0&0
\\ \noalign{\medskip}0&0&0&0&0&0&0&0&0&0&0&0&0&0&1&0&0&0&0&0&0&0
\\ \noalign{\medskip}0&0&0&0&0&0&0&0&0&0&0&0&0&0&0&0&0&0&0&1&0&0
\\ \noalign{\medskip}0&0&0&0&0&0&0&0&0&0&0&0&0&0&0&0&0&0&1&0&0&0
\\ \noalign{\medskip}0&0&0&0&0&0&0&0&0&0&0&0&0&0&0&0&0&1&0&0&0&0
\\ \noalign{\medskip}0&0&0&0&0&0&0&0&0&0&0&0&0&0&0&0&1&0&0&0&0&0
\\ \noalign{\medskip}0&0&0&0&0&0&0&0&0&0&0&0&0&0&0&1&0&0&0&0&0&0
\\ \noalign{\medskip}2&2&0&0&0&-1&-5&-4&-7&-10&-8&-6&-4&-2&-5&-4&-7&-
10&-8&-6&-4&-2\\ \noalign{\medskip}0&0&0&0&0&0&0&1&0&0&0&0&0&0&0&0&0&0
&0&0&0&0\end {array} \right] 
$$
}
\caption{$M(h_2)$}\label{fig:Mh2}
\end{figure}
\begin{figure}
{\tiny 
\setlength{\arraycolsep}{3.5pt} 
$$
\left[ \begin {array}{cccccccccccccccccccccc} 14&27&-9&9&9&-42&-30&-
24&-42&-60&-48&-36&-24&-12&0&0&0&0&0&0&0&0\\ \noalign{\medskip}8&14&-5
&5&5&-23&-15&-12&-21&-30&-24&-18&-12&-6&0&0&0&0&0&0&0&0
\\ \noalign{\medskip}5&9&-4&3&3&-14&-10&-8&-14&-20&-16&-12&-8&-4&0&0&0
&0&0&0&0&0\\ \noalign{\medskip}21&39&-13&12&13&-60&-40&-32&-56&-80&-64
&-48&-32&-16&0&0&0&0&0&0&0&0\\ \noalign{\medskip}-49&-84&28&-28&-29&
133&105&84&147&210&168&126&84&42&0&0&0&0&0&0&0&0\\ \noalign{\medskip}1
&3&-1&1&1&-5&0&0&0&0&0&0&0&0&0&0&0&0&0&0&0&0\\ \noalign{\medskip}0&0&0
&0&0&0&0&0&1&0&0&0&0&0&0&0&0&0&0&0&0&0\\ \noalign{\medskip}3&6&-2&2&2&
-9&-8&-5&-10&-15&-12&-9&-6&-3&0&0&0&0&0&0&0&0\\ \noalign{\medskip}0&0&0
&0&0&0&1&0&0&0&0&0&0&0&0&0&0&0&0&0&0&0\\ \noalign{\medskip}0&0&0&0&0&0
&0&0&0&1&0&0&0&0&0&0&0&0&0&0&0&0\\ \noalign{\medskip}0&0&0&0&0&0&0&0&0
&0&1&0&0&0&0&0&0&0&0&0&0&0\\ \noalign{\medskip}0&0&0&0&0&0&0&0&0&0&0&1
&0&0&0&0&0&0&0&0&0&0\\ \noalign{\medskip}0&0&0&0&0&0&0&0&0&0&0&0&1&0&0
&0&0&0&0&0&0&0\\ \noalign{\medskip}0&0&0&0&0&0&0&0&0&0&0&0&0&1&0&0&0&0
&0&0&0&0\\ \noalign{\medskip}0&0&0&0&0&0&0&0&0&0&0&0&0&0&1&0&0&0&0&0&0
&0\\ \noalign{\medskip}0&0&0&0&0&0&0&0&0&0&0&0&0&0&0&1&0&0&0&0&0&0
\\ \noalign{\medskip}0&0&0&0&0&0&0&0&0&0&0&0&0&0&0&0&1&0&0&0&0&0
\\ \noalign{\medskip}0&0&0&0&0&0&0&0&0&0&0&0&0&0&0&0&0&1&0&0&0&0
\\ \noalign{\medskip}0&0&0&0&0&0&0&0&0&0&0&0&0&0&0&0&0&0&1&0&0&0
\\ \noalign{\medskip}0&0&0&0&0&0&0&0&0&0&0&0&0&0&0&0&0&0&0&1&0&0
\\ \noalign{\medskip}0&0&0&0&0&0&0&0&0&0&0&0&0&0&0&0&0&0&0&0&1&0
\\ \noalign{\medskip}0&0&0&0&0&0&0&0&0&0&0&0&0&0&0&0&0&0&0&0&0&1
\end {array} \right] 
$$
}
\caption{$M(h_3)$}\label{fig:Mh3}
\end{figure}
The matrix representations $M(h_i)$ of $\tau(h_i)$ on $S_{X(7)}$ are given in
Figures~\ref{fig:Mh1}--\ref{fig:Mh3}. (Recall that $\OG(S_X)$ acts on $S_X$ from the right.
Hence $M(h_i)$ satisfies $M(h_i) \cdot G_{\Lambda} \cdot {}^t M(h_i)=G_{\Lambda}$,
where $G_{\Lambda}$ is the Gram matrix of $\Lambda_{7,1}^{-}$ with respect to $\ve_1, \dots, \ve_{22}$.)
The characteristic polynomial
of the product 
$$
M:=M(h_1)M(h_2)M(h_3)
$$
is a Salem polynomial  
\begin{eqnarray*}
&&t^{22}-993\,t^{21}-1152\,t^{20}-123\,t^{19}+924\,t^{18}+584\,t^{17}-500\,t^{16}-1022\,t^{15}\\
&&-661\,t^{14}+105\,t^{13}+476\,t^{12}+878\,t^{11}+476\,t^{10}+105\,t^9-661\,t^8\\
&&-1022\,t^7-500\,t^6+584\,t^5+924\,t^4-123\,t^3-1152\,t^2-993\,t+1, 
\end{eqnarray*}
which has a positive real root
$994.15889\dots$.
\bibliographystyle{plain}
\
\def\cftil#1{\ifmmode\setbox7\hbox{$\accent"5E#1$}\else
  \setbox7\hbox{\accent"5E#1}\penalty 10000\relax\fi\raise 1\ht7
  \hbox{\lower1.15ex\hbox to 1\wd7{\hss\accent"7E\hss}}\penalty 10000
  \hskip-1\wd7\penalty 10000\box7} \def\cprime{$'$} \def\cprime{$'$}
  \def\cprime{$'$} \def\cprime{$'$}

\end{document}